\documentclass{amsart}
\usepackage{amsmath} \usepackage{amssymb} \usepackage{euscript}
\usepackage[all]{xy} \usepackage{varioref}
\usepackage{verbatim}
\usepackage{xspace}
\theoremstyle{plain}
\newtheorem{theorem}{Theorem}

\theoremstyle{remark}
\newtheorem{definition}{Definition} \newtheorem{example}{Example}
\newtheorem{remark}{Remark}

\newcommand{\pd}[2]{\ensuremath{\frac{\partial{#1}}{\partial{#2}}}}
\newcommand{\R}[1]{\ensuremath{\mathbb{R}^{#1}}}

\newcommand{\SO}[1]{\operatorname{SO}\left({#1}\right)}

\begin{document}
\title{Sussmann's orbit theorem and maps} \author{Benjamin McKay}
\address{University College Cork \\
  Cork, Ireland} \email{B.McKay@UCC.ie}
\date{\today} 
\thanks{This work was supported in full or in part by a grant from 
the University of South Florida St. Petersburg New Investigator
Research Grant Fund. This support does not necessarily imply
endorsement by the University of research conclusions.}
\begin{abstract}
  A map between manifolds which matches up families of complete vector
  fields is a fiber bundle mapping on each orbit of those vector
  fields.
\end{abstract}
\maketitle
\section{Introduction}
\begin{definition}
  Write $e^{tX}(m) \in M$ for the flow of a vector field $X$ through a
  point $m$ after time $t$.  Let $\mathfrak{F}$ be a family of smooth
  vector fields on a manifold $M$.  The \emph{orbit} of $\mathfrak{F}$
  through a point $m \in M$ is the set of all points $e^{t_1X_1}
  e^{t_2X_2} \dots e^{t_k X_k}(m)$ for any vector fields $X_j \in
  \mathfrak{F}$ and numbers $t_j$ (positive or negative) for which
  this is defined.
\end{definition}
\begin{example}
  The vector field $\pd{}{\theta}$ on the Euclidean plane (in polar
  coordinates) has orbits the circles around the origin, and the
  origin itself.
\end{example}
\begin{example}
  The set of smooth vector fields supported in a disk has as orbits
  the open disk (a 2-dimensional orbit) and the individual points
  outside or on the boundary of the disk (zero dimensional orbits).
\end{example}
\begin{example}
  On Euclidean space, the set of vector fields supported inside a
  ball, together with the radial vector field coming from the center
  of the ball, forms a set of vector fields with a single orbit.
\end{example}
\begin{example}
  Translation in a generic direction on a flat torus has densely
  winding orbits.
\end{example}
H{\'e}ctor Sussmann \cite{Sussmann:1973a,Sussmann:1973b,Sussmann:1974}
proved that the orbits of any family of smooth vector fields are
immersed submanifolds.  We prove that a mapping between two manifolds
which carries one family of complete vector fields into another, is a
fiber bundle mapping on each orbit.
\section{Proofs}
For completeness, we prove Sussmann's theorem.
\begin{theorem}[Sussmann \cite{Sussmann:1973a}]\label{thm:orbit}
  The orbit of any point under any family of smooth vector fields is
  an immersed submanifold (in a canonical topology). If two orbits
  intersect, then they are equal.  Let $\bar{\mathfrak{F}}$ be the
  largest family of smooth vector fields which have the same orbits as
  the given family $\mathfrak{F}$. Then $\bar{\mathfrak{F}}$ is a Lie
  algebra of vector fields, and a module over the algebra of smooth
  functions.
\end{theorem}
\begin{remark}
  Obviously, one could localize these results, replacing globally
  defined vector fields with subsheaves of the sheaf of locally
  defined smooth vector fields.
\end{remark}
\begin{proof}
  We can replace $\mathfrak{F}$ by $\bar{\mathfrak{F}}$ without loss
  of generality. Therefore, if $X,Y \in \mathfrak{F}$, we can suppose
  that $e^X_* Y \in \mathfrak{F}$ since the flow of $e^X_*Y$ is \[
  e^{t \left(e^X_*Y\right)} = e^X e^{tY},
\]
which must preserve orbits. We refer to this process as \emph{pushing
  around} vector fields.

Fix attention on a specific orbit.  For each point $m_0 \in M$, take
as many vector fields as possible $X_1,\dots,X_k$, out of
$\mathfrak{F}$, which are linearly independent at $m$.  Refer to the
number $k$ of vector fields as the \emph{orbit dimension}.  Pushing
around convinces us that the orbit dimension is a constant throughout
the orbit.  Refer to the map
\[
\left(t_1, \dots, t_k\right) \in \text{open } \subset \R{k} \mapsto
e^{t_1 X_1} \dots e^{t_k X_k} m_0 \in M
\]
(which we will take to be defined in some open set on which it is an
embedding) as a \emph{distinguished chart} and its image as a
\emph{distinguished set}.  The tangent space to each point $e^{t_1
  X_1} \dots e^{t_k X_k} m_0$ of a distinguished set is spanned by the
linearly independent vector fields
\[
X_1,e^{t_1 X_1}_* X_2, \dots,e^{t_1 X_1} \dots e^{t_{k-1} X_{k-1}}_*
X_k,
\] 
which belong to $\mathfrak{F}$, since they are just pushed around
copies of the $X_j$.  Let $\Omega$ be a distinguished set. Suppose
that $Y \in \mathfrak{F}$ is a vector field, which is not tangent to
$\Omega$. Then at some point of $\Omega$, $Y$ is not a multiple of
those pushed around vector fields, so the orbit dimension must exceed
$k$.

Therefore all vector fields in $\mathfrak{F}$ are tangent to all
distinguished sets. So any point inside any distinguished set stays
inside that set under the flow of any vector field in $\mathfrak{F}$,
at least for a short time. So such a point must also stay inside the
distinguished set under compositions of flows of the vector fields, at
least for short time.  Therefore a point belonging to two
distinguished sets must remain in both of them under the flows that
draw out either of them, at least for short times. Therefore that
point belongs to a smaller distinguished set lying inside both of
them. Therefore the intersection of distinguished sets is a
distinguished set.

We define an open set of an orbit to be any union of distinguished
sets; so the orbit is locally homeomorphic to Euclidean space.  We can
pick a countable collection of distinguished sets as a basis for the
topology.  Every open subset of $M$ intersects every distinguished set
in a distinguished set, so intersects every open set of the orbit in
an open set of the orbit. Thus the inclusion mapping of the orbit into
$M$ is continuous.  Since $M$ is metrizable, the orbit is also
metrizable, so a submanifold of $M$.  The distinguished charts give
the orbit a smooth structure. They are smoothly mapped into $M$,
ensuring that the inclusion is a smooth map.
\end{proof}
\begin{example}
  Let $\alpha=dy-z \, dx$ in $\R{3}$.  The vector fields on which
  $\alpha=0$ have one orbit: all of $\R{3}$, since they include
  $\partial_z,\partial_x+z\partial_y$, and therefore include the
  bracket:
\[
\left[ \partial_z,\partial_x+z\partial_y \right] = \partial_y.
\]
\end{example}
\begin{definition}
  Take a map $\phi : M_0 \to M_1$, and vector fields $X_j$ on $M_j$,
  $j=0,1$.  Write $\phi_* X_0 = X_1$ to mean that for all $m_0 \in
  M_0$, $\phi'\left(m_0\right) X_0\left(m_0\right) =
  X_1\left(\phi\left(m_0\right)\right).$ For families of vector
  fields, write $\phi_* \mathfrak{F}_0 = \mathfrak{F}_1$ to mean that
\begin{enumerate}
\item for any $X_0 \in \mathfrak{F}_0$ there is an $X_1 \in
  \mathfrak{F}_1$ so that $\phi_* X_0 = X_1$ and
\item for any $X_1 \in \mathfrak{F}_1$ there is a vector field $X_0
  \in \mathfrak{F}_0$ so that $\phi_* X_0 = X_1$.
\end{enumerate}
\end{definition}
\begin{example}
  The vector field $\partial_x$ on $\R{}$ has $\R{}$ as orbit.
  Consider the inclusion $(0,1) \subset \R{}$ of some open interval.
  The orbit of $\partial_x$ on $(0,1)$ is $(0,1)$. The orbits are
  mapped to each other by the inclusion, but not surjectively.
\end{example}
\begin{example}
  If $M_0=\R{2}_{x,y}$ and $M_1=\R{1}_x$, and $\phi(x,y)=x$, and
  $\mathfrak{F}_0 = \left\{\partial_x,\partial_y\right\}$ and
  $\mathfrak{F}_1=\left\{\partial_x,0\right\},$ then clearly $\phi_*
  \mathfrak{F}_0 = \mathfrak{F}_1$.
\end{example}
\begin{example}
  The group $\SO{3}$ of rotations acts on the sphere $S^2$, and we can
  map $\SO{3} \to S^2$, taking a rotation $g$ to $gn$ where $n$ is the
  north pole. This map takes the left invariant vector fields to the
  infinitesimal rotations, and clearly is a fiber bundle, the Hopf
  fibration.
\end{example}
\begin{theorem}\label{cor:orbitMap}
  If $\mathfrak{F}_j$ are sets of vector fields on manifolds $M_j$,
  for $j=0,1$, and $\phi : M_0 \to M_1$ satisfies $\phi_*
  \mathfrak{F}_0 = \mathfrak{F}_1$, then $\phi$ takes
  $\mathfrak{F}_0$-orbits into $\mathfrak{F}_1$-orbits. On each orbit,
  $\phi$ has constant rank.  If the vector fields in both families are
  complete, then $\phi$ is a fiber bundle mapping on each orbit.
\end{theorem}
\begin{proof}
  By restricting to an orbit in $M_0$, we may assume that there is
  only one orbit. The map $\phi$ is invariant under the flows of the
  vector fields, so must have constant rank.
  
  Henceforth, suppose that the vector fields are complete.  Given a
  path
\[
e^{t_1 X_1} \dots e^{t_k X_k}m_0
\]
down in $M_1$, we can always lift it to one in $M_0$, so $\phi$ is
onto.  It might not be true that $\phi_* \bar{\mathfrak{F}}_0 =
\bar{\mathfrak{F}}_1$, but nonetheless we can still push around vector
fields, because the pushing upstairs in $M_0$ corresponds to pushing
downstairs in $M_1$. So without loss of generality, both
$\mathfrak{F}_0$ and $\mathfrak{F}_1$ are closed under ``pushing
around''.

As in the above proof, for each point $m_1 \in M_1$, we can construct
a distinguished chart
\[
\left(t_1, \dots, t_k\right) \mapsto e^{t_1 X_1} \dots e^{t_k X_k}
m_1.
\]
These $X_k$ are vector fields on $M_1$.  Write $Y_k$ for some vector
fields on $M_0$ which satisfy $\phi_* Y_k = X_k$.  Clearly $\phi$ is a
surjective submersion.  Let $U_1 \subset M_1$ be the associated
distinguished set; on $U_1$ these $t_j$ are now coordinates.  Let $U_0
= \phi^{-1} U_1 \subset M_0$. Let $Z$ be the fiber of $\phi : M_0 \to
M_1$ above the origin of the distinguished chart. Map
\[
u_0 \in U_0 \mapsto \left(u_1,z\right) \in U_1 \times Y
\]
by $u_1=\phi\left(u_0\right)$ and
\[
z = e^{-t_k Y_k} \dots e^{-t_1 Y_1} u_0.
\]
Clearly this gives $M_0$ the local structure of a product.  The
transition maps have a similar form, composing various flows, so $M_0
\to M_1$ is a fiber bundle.
\end{proof}
Keep in mind that all vector fields on compact manifolds are complete.
Even though the orbits might not be compact, our theorem says that the
orbits upstairs will fiber over the orbits downstairs.
\begin{example}
  Take $M_1 = \R{2}_{x,y}$, and $M_0 \subset M_1$ a pair of disjoint
  disks, say those of unit radius around two points of the $x$ axis
  which are distantly separated.  As the family $\mathfrak{F}_0$ up in
  the disks, take the translation vector fields
  $\partial_x,\partial_y$ along coordinate axes in the right disk, and
  in the left, the pair of vector fields $\partial_x,0$. Obviously
  these are not complete. As the family $\mathfrak{F}_1$, take the
  translation vector field $\partial_x$, and a vector field $f(x,y)
  \partial_y$ which vanishes in the left disk, and nowhere outside of
  closure of the left disk, and equals $\partial_y$ in the right disk.
  The orbits downstairs are all two dimensional, while those upstairs
  are one dimensional in the left disk, and two dimensional in the
  right.
\end{example}
\begin{example}
  Take $E \to M$ any fiber bundle, and pick a plane field on $E$
  transverse to the fibers. Every vector field on $M$ lifts to a
  unique vector field on $E$ tangent to the 2-plane field. Suppose
  that the fibers of $E \to M$ are compact. Lifting all complete
  vector fields, we get a family of complete vector fields on $E$.
  Their orbits must be connected and fiber over $M$.
\end{example}
\begin{example}
  Take any 2-plane field on $\SO{3}$ transverse to the leaves of the
  Hopf fibration $\SO{3} \to S^2$, and lift vector fields as in the
  last example.  A two dimensional orbit would have to be
  diffeomorphic to $S^2$, since $S^2$ is simply connected.  The Hopf
  fibration admits no section, so therefore all orbits must be three
  dimensional, hence open and disjoint, and cover $\SO{3}$, which is
  connected. Hence every 2-plane field transverse to the Hopf
  fibration has all of $\SO{3}$ as orbit, even though the 2-plane
  field may be holonomic on an open set.  The same result works for
  any circle bundle on any compact manifold: either every plane field
  transverse to the circle fibers has a single orbit, or the circle
  bundle trivializes on a covering space.
\end{example}
\begin{example}
  Consider the Hopf fibration $S^3 \to S^7 \to S^4$.  Take any 4-plane
  field on $S^7$ transverse to the fibers. The orbits must be bundles
  over $S^4$. The fibers of such a bundle cannot be zero dimensional,
  since $S^4$ is simply connected and the Hopf fibration is not a
  trivial bundle. Suppose that $F \to B \to S^4$ is a fiber bundle,
  and that $B \subset S^7$ is a subbundle.  The bundle $B$ cannot be
  trivial, since that would give rise to a section of the Hopf
  fibration. The bundle $B$ is determined completely by slicing $S^4$
  along the equatorial $S^3$, and mapping $S^3$ to the diffeomorphism
  group of the fiber $F$.  The fiber $F$ cannot be the real line, the
  circle, or a closed surface other than the sphere, since the
  diffeomorphism groups of these manifolds retract to finite
  dimensional groups which are aspherical. Therefore $F$ must be a
  sphere or noncompact surface, or a component of the complement in
  $S^3$ of a set of disjoint spheres and noncompact surfaces. Our
  theorem does not suffice to give a complete analysis of the possible
  orbits, but clearly it makes a substantial contribution to this
  question.
\end{example}
\bibliographystyle{amsplain} \bibliography{orbit}

\end{document}